\documentclass[twocolumn]{autart}    % Enable this line and disable the 
\usepackage{graphicx}          % Include this line if your 
\usepackage{cite}
\usepackage{psfrag}
\usepackage{amsmath}
\usepackage{amssymb}  % assumes amsmath package installed
\usepackage{xcolor}

\graphicspath{{Figures/}}

\newtheorem{defi}{Definition}

\newcommand{\tA}{\mathcal{A}}
\newcommand{\tB}{\mathcal{B}}
\newcommand{\tP}{\mathcal{P}}
\newcommand{\di}{\mathrm{diag}}
\newcommand{\tG}{\mathcal{G}}

\newcommand{\tL}{\mathcal{L}}
\newcommand{\tK}{\mathcal{K}}

\newcounter{MYtempeqncnt}

\begin{document}

\begin{frontmatter}

\title{Robust Consensus of Linear Multi-Agent Systems under Input Constraints or Uncertainties \thanksref{footnoteinfo}} 

\thanks[footnoteinfo]{This paper was not presented at any IFAC 
meeting. Corresponding author: D. H. Nguyen. }

\author{Dinh Hoa Nguyen}\ead{hoadn.ac@gmail.com}, %\ead{N.DinhHoa.Nguyen@ieee.org},    % Add the e-mail address 
\author{Tatsuo Narikiyo}\ead{n-tatsuo@toyota-ti.ac.jp},  % (ead) as shown
\author{Michihiro Kawanishi}\ead{kawa@toyota-ti.ac.jp},  % (ead) as shown

\address{Control System Laboratory, Department of Advanced Science and Technology, Toyota Technological Institute, 
        2-12-1 Hisakata, Tempaku-ku, Nagoya 468-8511, Japan.}% Please supply                                              
            % full addresses

%-----------------------------------------------------------------------------------          
                                          
%%%%%%%%%%%%%%%%%%%%%%%%%%%%%%%%%%%%%%%%%%%%%%%%%%%%%%%%%%%%%%%%%%%%%%%%%%%%%%%%
\begin{abstract}
This paper proposes a new approach to analyze and synthesize robust consensus control laws for general linear leaderless multi-agent systems (MASs) subjected to input constraints or uncertainties. 
First, the MAS under input constraints or uncertainties is reformulated as a network of Lur'e systems.    
Next, two scenarios of communication topology are considered, namely undirected and directed cyclic structures. 
In each case, a sufficient condition for consensus and the design of consensus controller gain are derived from solutions of a distributed LMI convex problem.
Finally, a numerical example is introduced to illustrate the effectiveness of the proposed theoretical approach. 
\end{abstract}

%-----------------------------------------------------------------------------------
\end{frontmatter}

%%%%%%%%%%%%%%%%%%%%%%%%%%%%%%%%%%%%%%%%%%%%%%%%%%%%%%%%%%%%%%%%%%%%%%%%%%%%%%%%
\section{Introduction}

Multi-agent systems (MASs) and their cooperative control problems have gained much attention since there are a lot of practical applications, e.g.,     
power grids, wireless sensor networks, transportation networks, systems biology, etc., can be formulated, analyzed and synthesized under the framework of MASs. 
A key feature in MASs is the achievement of a global objective by performing local measurement and control at each agent and simultaneously collaborating among agents using that local information.
One of the most important and intensively investigated issues in MASs (and their applications) is the consensus problem due to its attraction in both theoretical and applied aspects 
\cite{Olfati-Saber:2004,Olfati-Saber:2007,Ren:2007}.

Since all real control systems are subjected to physical constraints on their inputs or states, MASs are not exception. 
Therefore, the MAS consensus under input or state constraints is a significant and realistic problem. 
Likewise, the MAS consensus in presence of input or state uncertainties is realistic and worth studying. 
Practical examples include consensus of vehicles with limited speeds and limited working space, smart buildings energy control with temperature and humidity are required in specific ranges, 
just to name a few. 
However, most of the early researches on MASs were not aware of those practical issues, and 
%only focus on several practical issues on the communication among agents, e.g. switching interconnection protocol or time delays on communication channels. 
it has not been until recently that some studies have considered the cooperative control of MASs in presence of input or state constraints on each agent 
\cite{Nedic:2010,Lin:2014,Lee:2011,Liu:2012,Wang:2013,Meng:2012,Su:2013,Takaba:2014,Zhang:2014,Zhang:2014j}.  

A constrained consensus problem was investigated in \cite{Nedic:2010} where the states of agents are required to lie in individual closed convex sets and the final consensus state must belong to the non-empty intersection of those sets. Accordingly, a projected consensus algorithm was proposed and then applied to distributed optimization problems. Following this research line, \cite{Lin:2014} extended the result in \cite{Nedic:2010} to the context where communication delays exist.  
In another work, \cite{Lee:2011} studied the state increment by utilizing the model predictive control (MPC) method. In fact, using the MPC framework we can also incorporate input or state constraints,  however the computational cost could be high. Therefore, distributed and fast MPC algorithms need to be developed to fit into the context of large-scale MASs. 
Another direction to deal with input or state constraints is to employ the so-called discarded consensus algorithms \cite{Liu:2012,Wang:2013}. Nevertheless, a disadvantage of these approaches as well as in \cite{Nedic:2010,Lin:2014} is that the initial states of agents must belong to some sets specified by the constraints, or in other words the consensus is only local. Moreover, only agents with single integrator dynamics were considered in \cite{Liu:2012,Wang:2013}.  

To achieve the global or semi-global consensus in presence of input or state constraints, some consensus laws were presented in \cite{Meng:2012,Su:2013}, but they were only for leader-follower MASs. Another way to tackle the input or state constraints to derive global consensus is to reformulate the constrained MAS as a network of Lur'e systems 
\cite{Takaba:2014,Zhang:2014,Zhang:2014j}. 
The paper \cite{Takaba:2014} considered linear agents with bounded-constraint inputs and obtained a sufficient condition for global consensus, but agents is limited to be single-input and the network 
%inter-agent communication 
is undirected. 
Next, \cite{Zhang:2014} and \cite{Zhang:2014j} investigated consensus problems where outputs of agents are incrementally bounded or passive, 
with directed and undirected topologies. Consequently, sufficient conditions for global consensus were derived in the form of LMI convex problems.  

Following the ideas of achieving global consensus by reformulating the considered MAS as a network of Lur'e systems in \cite{Takaba:2014,Zhang:2014,Zhang:2014j}, this paper proposes a new approach to design robust consensus controllers for general linear homogeneous leaderless MASs in presence of input constraints or uncertainties. 
%where the input vector of each agent is subjected to a box constraint.  
%By employing Lyapunov theory and Lasalle's invariance principle, sufficient conditions for robust consensus are derived in the form of LMI problems.   
The contributions of this paper are threefold. First, the proposed approach is applicable for leaderless MASs with general linear dynamics of agents, and the class of nonlinearities induced by 
the input constraints or uncertainties is broader than those in \cite{Takaba:2014,Zhang:2014,Zhang:2014j}. 
Second, the consensus controller gain is computed from the solution of a distributed low-dimension convex problem with LMI constraints. 
%which are more general and flexible then the one in \cite{Takaba:2014}. 
Third, the proposed approach can be used for global consensus analysis and synthesis under both scenarios of undirected networks and a special class of directed networks.
%input and state constraints since the latter case can be transformed into an equivalent former one.  

%The paper is organized as follows. In section \ref{prob-form}, we introduce a class of Goodwin oscillator networks 
%which is employed to describe the circadian oscillations. We then propose an approach based on harmonic balance 
%method to investigate the phases and amplitudes dependence of network oscillations on exogenous inputs in section 
%\ref{hb-method}. In section \ref{HO-osc}, a further study on Goodwin oscillator networks excited by periodic inputs 
%with higher order harmonics is presented. Lastly, our main results are summarized in section \ref{sum}.

The following notation and symbols will be used in the paper. $\mathbb{R}$ and $\mathbb{C}$ stand for the real 
and complex sets, and $j$ denotes the complex unit.  
%$\mathbb{R}^n$, $\mathbb{R}^n_{+}$ and $\mathbb{C}^n$ are used to denote the set of real, 
%positive real and complex $n\times 1$  vectors. 
Moreover, $\mathbf{1}_n$ denotes the $n\times1$ vector 
with all elements equal to $1$, and $I_n$ denotes the $n\times n$ identity matrix.
Next, $\otimes$ stands for the Kronecker product, $\di\{\}$ denotes diagonal or block-diagonal matrices, and   
$\mathrm{sym}(A)$ denotes $A+A^T$ for any real matrix $A$.
Lastly, $\succ$ and $\succeq$ denote the positive definiteness and positive semi-definiteness of a matrix, and similar meanings are used for $\prec$ and $\preceq$.

%============================================================================================================================
\section{Problem Setting}
\label{prob}

%----------------------------------------------------------------------
\subsection{Graph Theory}
%\label{graph}

Denote $(\mathcal{G},\mathcal{V},\mathcal{E})$ the graph representing the information structure in an MAS composing of $N$ agents, 
where each node in $\mathcal{G}$ stands for an agent and each edge in $\mathcal{G}$ represents the interconnection between two agents; 
$\mathcal{V}$ and $\mathcal{E}$ represent the set of vertices and edges of $\mathcal{G}$, respectively. There is an edge $e_{ij} \in \mathcal{E}$ if agent $i$ receives information from agent $j$. 
The neighboring set of a vertex $i$ is denoted by $\mathcal{N}_{i} \triangleq \{j:e_{ij} \in \mathcal{E}\}$.  
Moreover, let $a_{ij}$ be elements of the adjacency matrix $\mathcal{A}$ of $\mathcal{G}$, i.e., $a_{ij}>0$ if $e_{ij} \in \mathcal{E}$ and  $a_{ij}=0$ if $e_{ij} \notin \mathcal{E}$. 
The in-degree of a vertex $i$ is denoted by $\mathrm{deg}^{\rm in}_{i} \triangleq \sum_{j=1}^{N}{a_{ij}}$, then the in-degree matrix of $\mathcal{G}$ is denoted by 
$\mathcal{D}=\mathrm{diag}\{\mathrm{deg}^{\rm in}_{i}\}_{i=1,\ldots,N}$. 
Consequently, the Laplacian matrix $\mathcal{L}$ associated to $\mathcal{G}$ is defined by $\mathcal{L}=\mathcal{D}-\mathcal{A}$. 
The out-degree of a vertex $i$ is denoted by $\mathrm{deg}^{\rm out}_{i} \triangleq \sum_{j=1}^{N}{a_{ji}}$. 
Then $\mathcal{G}$ is said to be balanced if $\mathrm{deg}^{\rm in}_{i}=\mathrm{deg}^{\rm out}_{i} \;\forall\; i=1,\ldots,N.$ 
A directed path connecting vertices $i$ and $j$ in $\mathcal{G}$ is a set of consecutive edges starting from $i$ and stopping at $j$. 
Then $\mathcal{G}$ is said to have a spanning tree if there exists a node called root node from which there are directed paths to every other node. 
$\tG$ is undirected if and only if $a_{ij}=a_{ji} \, \forall \, i,j=1,\ldots,N$.

%----------------------------------------------------------------------
\subsection{Problem Description}
%\label{}

Consider a MAS including of $N$ identical agents with the following linear dynamics
\begin{equation}
	\label{agent}
		\dot{x}_i = Ax_i+Bu_i,  i=1,\ldots,N,
\end{equation} 
where $x_i \in \mathbb{R}^{n}$ is the state vector, $u_i \in \mathbb{R}^{m}$ is the control input,  
$A \in \mathbb{R}^{n\times n}$, $B \in \mathbb{R}^{n \times m}$. 
The following assumptions will be employed. 
\begin{itemize}
	\item[\bf A1:]  $(A,B)$ is stabilizable. 
	\item[\bf A2:]  All eigenvalues of $A$ is on the closed left half complex plane.
	\item[\bf A3:]  $\tG$ is balanced and contains a spanning tree.
\end{itemize}
Assumptions A1-A2 are necessary and sufficient such that the consensus can be achieved and stable (see e.g. \cite{Ma:2010}). 
And assumption A3 implies that $\tG$ is connected if it is undirected. 
 
Denote $x = \left[x_{1}^T,\ldots,x_{N}^T\right]^T, u=\left[u_{1}^T,\ldots,u_{N}^T\right]^T$. 
The whole MAS at the initial state is then described by
\begin{equation}
	\label{mas}
		\dot{x} = \tA x+\tB u.
\end{equation}
It has recently been proved in our previous research \cite{Nguyen:2015TAC1} that without any further requirement on the control input or agents' states, the MAS (\ref{mas}) can reach consensus in the sense of (\ref{consensus-defi}) by a control law in the following form,
\begin{equation}
	\label{consensus-ctlr}
	u=(\tL \otimes \tK)x,
\end{equation}
with a properly synthesized $\tK$. 
%where $\tL$ is the Laplacian matrix associated with the information graph in the MAS, $\tK=-\mu RB^TP$ with any $\mu>0$ and $P\in\mathbb{R}^{n\times n},P \succ 0$ is the unique solution of the following Riccati equation
%\begin{equation}
	%\label{Req-local}
	%PA+A^TP+Q-PBRB^TP=0,
%\end{equation}
%of which $Q\in\mathbb{R}^{n\times n}$, $Q\succeq 0$ such that $(Q^{1/2},A)$ is detectable, $R\in\mathbb{R}^{m\times m}$, $R\succ 0$. 
%Nevertheless, in real applications, due to physical limitations of agents and the uncertain communication links among them, the states and inputs of agents are usually bounded in some certain ranges 
%and may contain some uncertainties. As a result, the control law (\ref{consensus-ctlr}) can no longer guarantee the consensus of agents. 
Nevertheless, in real applications the inputs of agents are usually bounded in some certain ranges due to physical limitations of agents, and may contain some uncertainties because of uncertain communication links.  
As a result, the control law (\ref{consensus-ctlr}) can no longer guarantee the consensus of agents. 
Therefore, to take into account the aforementioned practical issues, we will consider in this research the following control scenario:
\begin{itemize}
	\item {\bf Input constraint/uncertainty:} For all $i \in [1,N]$, $u_{i,k}=f_{i,k}(z_{i,k}) \,\forall\, k=1,\ldots,m$ where $z_{i,k} \in \mathbb{R}$ is the aggregated signal that the $k$th 
	input of agent $i$ received; $f_{i,k}: \mathbb{R} \rightarrow \mathbb{R}$ is a continuous function that satisfies the following sector-bounded condition:
		\begin{equation}
			\label{input-cstrt}
			\begin{aligned}
				&(f_{i,k}(z_{i,k})-\delta_{k,1}z_{i,k})(f_{i,k}(z_{i,k})-\delta_{k,2}z_{i,k}) \leq 0, \\ 
				& \, \forall \, k=1,\ldots,m; \, \forall \, i=1,\ldots,N,
			\end{aligned}
		\end{equation}
		where $\delta_{k,1},\delta_{k,2} \in \mathbb{R}$ are known constants, $\delta_{k,1}<\delta_{k,2}$.
	%\item[{\bf (ii)}] {\bf State constraint/uncertainty:} For all $i \in [1,N]$, $y_{i,k}=g_{i,k}(x_{i,k}) \,\forall\, k=1,\ldots,n$ where $y_{i,k} \in \mathbb{R}$ is the $k$th component of the signal sent 
	%from agent $i$; $g_{i,k}: \mathbb{R} \rightarrow \mathbb{R}$ is a continuous function that satisfies the following sector-bounded condition:
		%\begin{equation}
			%\label{state-cstrt}
			%\begin{aligned}
				%&(g_{i,k}(x_{i,k})-w_{k,1}x_{i,k})(g_{i,k}(x_{i,k})-w_{k,2}x_{i,k}) \leq 0 \\ 
				%& \, \forall \, k=1,\ldots,m; \, \forall \, i=1,\ldots,N,
			%\end{aligned}
		%\end{equation}
		%where $w_{k,1},w_{k,2} \in \mathbb{R}$ are known constants, $w_{k,1}<w_{k,2}$. 
\end{itemize}
Consequently, in presence of input constraints or uncertainties described above, each agent try to collaborate with others to achieve a consensus defined as follows. 
%Let us denote $\tG$ the graph representing the information structure in the MAS, 
%in which each node in $\mathcal{G}$ represents an agent and each edge in $\mathcal{G}$ represents 
%the interconnection between two agents. 
%%We assume that $\mathcal{G}$ is undirected and connected. 
%Consequently, the consensus of agents is defined as follows. 
\begin{defi}
The MAS with linear dynamics of agents represented by (\ref{agent}) and the information exchange among agents represented by 
$\tG$ is said to reach a consensus if 
\begin{equation}
	\label{consensus-defi}
	\lim_{t \rightarrow \infty} \|x_{i}(t)-x_{j}(t)\| = 0 ~\forall~i,j=1,\ldots,N.
\end{equation}
\end{defi}
Next, we introduce the control design problem investigated in this paper. 
%\begin{itemize}
	%\item {\bf Design problem 1 (Robust consensus under input constraint/uncertainty):} For the given linear MASs with dynamics of agents represented by (\ref{agent}) and the information exchange among agents represented by $\tG$, find a control strategy to achieve consensus of agents in the sense of (\ref{consensus-defi}) under the input constraint/uncertainty (\ref{input-cstrt}). 
	%\item {\bf Design problem 2 (Robust consensus under state constraint/uncertainty):} For the given linear MASs with dynamics of agents represented by (\ref{agent}) and the information exchange among agents represented by $\tG$, find a control strategy to achieve consensus of agents in the sense of (\ref{consensus-defi}) under the state constraint/uncertainty (\ref{state-cstrt}).
%\end{itemize}
\begin{itemize}
	\item {\bf Design problem  (Robust consensus under input constraints or uncertainties):} For the given linear MASs with dynamics of agents represented by (\ref{agent}) and the information exchange among agents represented by $\tG$, find a control strategy to achieve consensus of agents in the sense of (\ref{consensus-defi}) under the input constraints or uncertainties (\ref{input-cstrt}), 
	for any initial conditions of agents. 
\end{itemize}

%============================================================================================================================
\section{Robust Consensus Analysis and Design under Input Constraints or Uncertainties}
\label{input}

Under the input constraints or uncertainties (\ref{input-cstrt}), we propose to use the following consensus control law 
\begin{equation}
	\label{input-cstrt-consensus-ctlr}
	u=\mathbf{f}\left( (\tL \otimes \tK)x \right), 
\end{equation}
where $\mathbf{f}(y) \triangleq [f_{1}^T(y_{1}),\ldots,f_{N}^T(y_{N})]^T$, $f_{i}(y_{i}) \triangleq [f_{i,1}^T(y_{i,1})$, $\ldots,f_{i,N}^T(y_{i,N})]^T$, $\forall\; y=[y_{1}^T,\ldots,y_{N}^T]^T$.  
Then the MAS (\ref{mas}) with this control strategy can be rewritten in the following form
\begin{equation}
	\label{mas-rewrite}
	\begin{aligned}
		\dot{x} &= \tA x+\tB u, \\
		z &= (\tL \otimes \tK)x, \\
		u &= \mathbf{f}(z),
	\end{aligned}
\end{equation}
which can be seen as a network of Lur'e systems. Note that this Lur'e network is different from that in \cite{Zhang:2014,Zhang:2014j} and the nonlinearity is more general.  

The following theorem presents a sufficient condition for robust consensus under input constraints or uncertainties and how to design the consensus controller gain $\tK$. 

%%%%%%%%%%%%%%%%%%%%%%%%%%%%%%
\begin{thm}
\label{input-consensus-thm}
The robust consensus is achieved for the MAS (\ref{mas}) with an {\bf undirected} communication graph by the control law (\ref{input-cstrt-consensus-ctlr}) 
if there exist matrices $X \in \mathbb{R}^{n\times n}$, $Y \in \mathbb{R}^{m\times n}$ and $Z \in \mathbb{R}^{m\times m}$ 
such that the following LMI problem is feasible with a given $\epsilon > 0$,
\begin{equation}
	\label{input-LMI}
	\scalebox{0.9}{$
	\begin{aligned}
		& \begin{bmatrix} \mathrm{sym}(AX+\lambda_{2}B\Delta_{2}Y)+\epsilon X & BZ+\frac{1}{2}\lambda_{2}Y^T(\Delta_{1}-\Delta_{2}) \\ 
		\left(BZ+\frac{1}{2}\lambda_{2}Y^T(\Delta_{1}-\Delta_{2})\right)^T & -Z \end{bmatrix} \preceq 0, \\
		& \begin{bmatrix} \mathrm{sym}(AX+\lambda_{N}B\Delta_{2}Y)+\epsilon X & BZ+\frac{1}{2}\lambda_{N}Y^T(\Delta_{1}-\Delta_{2}) \\ 
		\left(BZ+\frac{1}{2}\lambda_{N}Y^T(\Delta_{1}-\Delta_{2})\right)^T & -Z \end{bmatrix} \preceq 0, \\
 		& X \succ 0, \\
		& Z \succ 0, Z \;\textrm{is diagonal},		
	\end{aligned}  $}
\end{equation}
where $\Delta_{1}=\di\{\delta_{k,1}\}_{k=1,\ldots,m}$, $\Delta_{2}=\di\{\delta_{k,2}\}_{k=1,\ldots,m}$, 
$\lambda_{2} \leq \lambda_{3} \leq \cdots \leq \lambda_{N}$ are non-zero eigenvalues of $\tL$. 
Furthermore, the controller gain $\tK$ is calculated by
\begin{equation}
	\label{gain}
	\tK=YX^{-1}.
\end{equation}
\end{thm}
%%%%%%%%%%%%%%%%%%%%%%%%%%%%%%

%%%%%%%%%%%%%%%%%%%
\begin{figure*}[!t]
% ensure that we have normalsize text
\normalsize
% Store the current equation number.
\setcounter{MYtempeqncnt}{\value{equation}}
% Set the equation number to one less than the one
% desired for the first equation here.
% The value here will have to changed if equations
% are added or removed prior to the place these
% equations are referenced in the main text.
\setcounter{equation}{11}
\begin{align}
	\label{input-LMI-d}
	\begin{bmatrix} \mathrm{sym}(AX+\frac{1}{2}B\tilde{\Delta}_{i}Y)+\epsilon X & BZ+Y^T\hat{\Delta}_{i} & -\mathbb{Y}_{i} & 0 \\ Z^TB^T+\hat{\Delta}_{i}Y & -Z & 0 & 0 \\
	\mathbb{Y}_{i} & 0 & \mathrm{sym}(AX+\frac{1}{2}B\tilde{\Delta}_{i}Y)+\epsilon X & BZ+Y^T\hat{\Delta}_{i} \\ 0 & 0 & Z^TB^T+\hat{\Delta}_{i}Y & -Z 	\end{bmatrix} \preceq 0, \; \forall \; i=2,\ldots,N.
\end{align}

% Restore the current equation number.
\setcounter{equation}{\value{MYtempeqncnt}}
% IEEE uses as a separator
\hrulefill
% The spacer can be tweaked to stop underfull vboxes.
\vspace*{4pt}
\end{figure*}	
%%%%%%%%%%%%%%%%%%%

\begin{pf}
Consider a Lyapunov function $V(x)=x^T\tP x$ where $\tP \triangleq P_{1}\otimes P_{2}$, $P_{1} \in \mathbb{R}^{N}$, $P_{1} \succeq 0$, $P_{2} \in \mathbb{R}^{n}$, $P_{2} \succ 0$. 
Taking the derivative of $V(x)$ gives us
\begin{equation*}
	\dot{V}(x) = x^T\left( \tP\tA+\tA^T\tP \right)x+2x^T\tP\tB u. 
\end{equation*}
Hence, for all $\epsilon > 0$ we have
\begin{equation*}
		\dot{V}(x)+\epsilon V(x) = %\begin{bmatrix} x \\ u \end{bmatrix}^T \begin{bmatrix} \tP\tA+\tA^T\tP+\epsilon\tP & \tP\tB \\ \tB^T\tP & 0 \end{bmatrix} \begin{bmatrix} x \\ u \end{bmatrix}. 
		x^T\left( \tP\tA+\tA^T\tP+\epsilon\tP \right)x+2x^T\tP\tB u.
\end{equation*} 
We now seek $\tP$ such that $\dot{V}(x)+\epsilon V(x) \leq 0$ as long as (\ref{input-cstrt}) holds. 
Using the S-procedure \cite{Boyd:2004}, such $\tP$ exists if there exist $\gamma_{1,1},\ldots,\gamma_{1,m},\ldots,\gamma_{N,1},\ldots,\gamma_{N,m}$ which are non-negative such that 
\begin{align}
	\label{Vdot-1}
	& \dot{V}(x)+\epsilon V(x)-\sum_{i=1}^{N}\sum_{k=1}^{m}{\gamma_{i,k}(u_{i,k}-\delta_{k,1}z_{i,k})(u_{i,k}-\delta_{k,2}z_{i,k})} \nonumber \\
	&  \leq 0, 
\end{align}
Let $\gamma_{i,k}=\gamma_{k} > 0 \,\forall\, i=1,\ldots,N$ and $\Gamma=\mathrm{diag}\{\gamma_{k}\}_{k=1,\ldots,m}$, then (\ref{Vdot-1}) is satisfied if 
\begin{align}
	\label{Vdot-2}
	& \dot{V}(x)+\epsilon V(x) -\sum_{i=1}^{N}{(u_{i}-\Delta_{1}z_{i})^T\Gamma(u_{i}-\Delta_{2}z_{i})} \leq 0, \nonumber \\
	& \Leftrightarrow \begin{bmatrix} x \\ u \end{bmatrix}^T \begin{bmatrix} \mathbb{P}_{1} & \mathbb{P}_{2} \\ \mathbb{P}_{2}^T & \mathbb{P}_{3} \end{bmatrix} \begin{bmatrix} x \\ u \end{bmatrix} \preceq 0 
	\Leftrightarrow \begin{bmatrix} \mathbb{P}_{1} & \mathbb{P}_{2} \\ \mathbb{P}_{2}^T & \mathbb{P}_{3} \end{bmatrix} \preceq 0,
\end{align}
where $\mathbb{P}_{1}=\tP\tA+\tA^T\tP+\epsilon\tP-(\tL^T\tL)\otimes(\tK^T\Gamma \Delta_{1}\Delta_{2}\tK)$, $\mathbb{P}_{2}=\tP\tB+\frac{1}{2}\tL^T\otimes(\tK^T\Gamma(\Delta_{1}+\Delta_{2}))$, 
$\mathbb{P}_{3}=-I_{N}\otimes \Gamma$.

%%%%%%%%%%%%%%%%%%%
\addtocounter{equation}{1}
%%%%%%%%%%%%%%%%%%%

Subsequently, employing Schur complement \cite{Boyd:2004} to (\ref{Vdot-2}) results in 
$\mathbb{P}_{1}-\mathbb{P}_{2}\mathbb{P}_{3}^{-1}\mathbb{P}_{2}^T \preceq 0$, which is equivalent to
\begin{align}
	\label{Vdot-3}
	%& \mathbb{P}_{1}-\mathbb{P}_{2}\mathbb{P}_{3}^{-1}\mathbb{P}_{2}^T \preceq 0, \nonumber \\
	%\Leftrightarrow 
	& P_{1}\otimes(A^TP_{2}+P_{2}A+\epsilon P_{2})+ P_{1}^2\otimes(P_{2}B\Gamma^{-1}B^TP_{2}) \nonumber \\
	& +\mathrm{sym}(\frac{1}{2}(P_{1}\tL)\otimes[P_{2}B(\Delta_{1}+\Delta_{2})\tK]) \nonumber \\ %+\frac{1}{2}(\tL^TP_{1})\otimes[\tK^T(S_{1}+S_{2})B^TP_{2}] \nonumber \\
	& +\frac{1}{4}(\tL^T\tL)\otimes[\tK^T\Gamma(\Delta_{1}-\Delta_{2})^2\tK] \preceq 0.
\end{align}

Let us choose $P_{1}=I_{N}-\frac{1}{N}\mathbf{1}_{N}\mathbf{1}_{N}^T$ then we can easily show that $P_{1}^2=P_{1}$, $\tL^TP_{1}=\tL^T$, and $P_{1}\tL=\tL$ for balanced graphs. 
Therefore, (\ref{Vdot-3}) is equivalent to 
\begin{align}
	\label{Vdot-4}
	& P_{1}\otimes(A^TP_{2}+P_{2}A+\epsilon P_{2}+P_{2}B\Gamma^{-1}B^TP_{2}) \nonumber \\
	& +\frac{1}{4}(\tL^T\tL)\otimes[\tK^T\Gamma(\Delta_{1}-\Delta_{2})^2\tK] \nonumber \\
	& +\mathrm{sym}(\frac{1}{2}\tL\otimes[P_{2}B(\Delta_{1}+\Delta_{2})\tK])  \preceq 0. %+\frac{1}{2}\tL^T\otimes[\tK^T(S_{1}+S_{2})B^TP_{2}]
\end{align}
Next, denote $X \triangleq P_{2}^{-1}$ and multiply $I_{N}\otimes X$ both to the left and to the right of (\ref{Vdot-4}), we obtain
\begin{align}
	\label{Vdot-5}
	& P_{1}\otimes(XA^T+AX+\epsilon X+B\Gamma^{-1}B^T) \nonumber \\
	& +\frac{1}{4}(\tL^T\tL)\otimes[Y^T\Gamma(\Delta_{1}-\Delta_{2})^2Y] \nonumber \\
	& +\mathrm{sym}(\frac{1}{2}\tL\otimes[B(\Delta_{1}+\Delta_{2})Y]) \preceq 0, %+\frac{1}{2}\tL^T\otimes[Y^T(S_{1}+S_{2})B^T] \preceq 0,
\end{align}
where $Y \triangleq \tK X$. 
For undirected graph $\tG$, $\tL=\tL^T$, so let us denote $U\in\mathbb{R}^{N\times N}$ the orthogonal matrix that diagonalizes $\tL$.  
Accordingly, applying a congruence transformation with $U\otimes I_{n}$ to (\ref{Vdot-5}) gives us
\begin{align}
	\label{Vdot-6}
	& XA^T+AX+\epsilon X+B\Gamma^{-1}B^T+\frac{1}{4}\lambda_{i}^2Y^T\Gamma(\Delta_{1}-\Delta_{2})^2Y \nonumber \\
	& +\mathrm{sym}(\frac{1}{2}\lambda_{i}B(\Delta_{1}+\Delta_{2})Y) \preceq 0, %+\frac{1}{2}\lambda_{i}Y^T(S_{1}+S_{2})B^T \preceq 0,
\end{align}
for all $i=2,\ldots,N$, since $U^TP_{1}U=\di\{0,1,\ldots,1\}$, $U^T\tL U=\di\{0,\lambda_{2},\ldots,\lambda_{N}\}$, and $U^T\tL^2 U=\di\{0,\lambda_{2}^2,$ 
$\ldots,\lambda_{N}^2\}$. 
By some simple mathematical manipulations, we can rewrite (\ref{Vdot-6}) as follows,
\begin{align}
	\label{Vdot-7}
	& XA^T+AX+\epsilon X+\lambda_{i}(Y^T\Delta_{2}B^T+B\Delta_{2}Y) \nonumber \\
	& +(BZ+\frac{1}{2}\lambda_{i}Y^T(\Delta_{1}-\Delta_{2}))Z^{-1}(ZB^T+\frac{1}{2}\lambda_{i}(\Delta_{1}-\Delta_{2})Y) \nonumber \\
	& \preceq 0,
\end{align} 
where $Z \triangleq \Gamma^{-1}$. 
Then using Schur complement again with (\ref{Vdot-7}) and noting that $\lambda_{i},i=3,\ldots,N-1$ can be represented as convex combinations of $\lambda_{2}$ and $\lambda_{N}$ since 
$\lambda_{2} \leq \lambda_{3},\ldots,\lambda_{N-1} \leq \lambda_{N}$, we obtain (\ref{input-LMI}). 
Next, we have seen that (\ref{Vdot-6}) implies (\ref{Vdot-1}) and hence 
\begin{equation}
	\label{inva-set}
	\begin{aligned}
		& \dot{V}(x)+\epsilon V(x) \\
		& \leq \sum_{i=1}^{N}\sum_{k=1}^{m}{\gamma_{i,k}(u_{i,k}-\delta_{k,1}z_{i,k})(u_{i,k}-\delta_{k,2}z_{i,k})} \leq 0, \\ 
		&\Rightarrow  \dot{V}(x) \leq -\epsilon V(x) \leq 0.	
	\end{aligned}
\end{equation}
Therefore, from Lasalle's invariance principle we can conclude that $\xi$ globally exponentially converges to the largest invariance set contained in 
$\left\{ \left. x \in \mathbb{R}^{nN} \right| \dot{V}(x)=0 \right\}$ for any initial condition. 
Furthermore, it can be seen from (\ref{inva-set}) that $\dot{V}(x)=0$ if and only if $V(x)=0$ which is equivalent to $x=\mathbf{1}_{N}\otimes \bar{x}$, $\bar{x} \in \mathbb{R}^{n}$, 
i.e., the consensus is achieved. 
\end{pf}

\begin{rem}
When the inputs of agents are subjected to boundedness, $f_{i}$ becomes the vector-valued saturation functions and hence $\delta_{k,1}=0$, $\delta_{k,2}=1$ $\forall \, k=1,\ldots,m$. 
%$w_{k,1}=0$, $w_{k,2}=1$ $\forall \, k=1,\ldots,n$. 
This particular case was investigated in \cite{Nguyen-ECC16} by a different control design. 
The method presented in this paper is more general and is applicable for more contexts than the one in \cite{Nguyen-ECC16}.    
\end{rem}

\begin{rem}
%The LMI problem (\ref{input-LMI}) can be solved in a distributed manner if $\lambda_{2}$ and $\lambda_{N}$ are known or can be approximated. 
Recently, there are several existing researches, e.g. \cite{Franceschelli:2013}, \cite{Tran:2015}, which propose different distributed methods to approximate the whole eigen-spectrum of the Laplacian matrix. These methods can be employed to estimate $\lambda_{2}$ and $\lambda_{N}$ before solving the LMI problem (\ref{input-LMI}).  
As a result, we can solve (\ref{input-LMI}) in a distributed fashion.   
\end{rem}

\begin{rem}
The results in \cite{Takaba:2014} can be considered as a special case of our result in Theorem \ref{input-consensus-thm} with single-input agents, input saturation,  and $Z=I_{m}$. 
Our result are much more general with the following properties: 
(i) its robustness to any constraint or uncertainty specified by (\ref{input-cstrt}); 
(ii) its applicability for leaderless MASs with general linear dynamics of agents; 
(iii) an additional variable $Z$ is introduced in the LMI problem (\ref{input-LMI}), which makes the LMIs less restrictive (cf. identity matrix in LMI problems (8) and (12) in \cite{Takaba:2014}); %, e.g., for single integrator dynamics (i.e., $A=0,B=1$), the LMI problems in \cite{Takaba:2014} (cf. (8) and (12) in \cite{Takaba:2014}) are infeasible but our LMIs result in $Y<0$ and any $X>0$; 
(iv) the term $\epsilon X$ in the upper-left blocks of the matrices in the LMI problem (\ref{input-LMI}) makes the consensus speed faster since $\dot{V}(x)$ is exponentially converged instead of being asymptotically converged as in \cite{Takaba:2014}.  
\end{rem}

Next, we present a design for directed networks in the following theorem. 

%%%%%%%%%%%%%%%%%%%%%%%%%%%%%%
\begin{thm}
\label{input-consensus-thm-1}
The robust consensus is achieved for the MAS (\ref{mas}) with a {\bf directed cyclic} unweighted communication graph by the control law (\ref{input-cstrt-consensus-ctlr}) 
if there exist matrices $X \in \mathbb{R}^{n\times n}$, $Y \in \mathbb{R}^{m\times n}$ and $Z \in \mathbb{R}^{m\times m}$ 
such that the LMI problem (\ref{input-LMI-d}) is feasible with a given $\epsilon > 0$, where for all $i=2,\ldots,N$,
\begin{equation}
	%\label{eq:}
	\begin{aligned}
		\mathbb{Y}_{i} &\triangleq -\frac{1}{2}\sin\frac{2\pi(i-1)}{N}\mathrm{sym}(B(\Delta_{1}+\Delta_{2})Y), \\
		\hat{\Delta}_{i} &\triangleq \sqrt{2\left(1-\cos\frac{2\pi(i-1)}{N}\right)}(\Delta_{1}-\Delta_{2}), \\
		\tilde{\Delta}_{i} &\triangleq \left(1-\cos\frac{2\pi(i-1)}{N}\right)(\Delta_{1}+\Delta_{2})-\hat{\Delta}_{i}.
	\end{aligned}	
\end{equation}
%and $\lambda_{i,r}$, $\lambda_{i,\ell}$ are the real and imaginary parts of the eigenvalues $\lambda_{i}$ of $\tL$.  
Accordingly, the controller gain $\tK$ is calculated by (\ref{gain}). 
\end{thm}
%%%%%%%%%%%%%%%%%%%%%%%%%%%%%%

\begin{pf}
Here, we employ the same Lyapunov function as in the proof of Theorem \ref{input-consensus-thm}, so all the steps until Eq. (\ref{Vdot-5}) are also applied for this scenario. 
Afterward, we note that $\tL$ is a circulant matrix since $\tG$ is an unweighted directed cyclic graph.  
Therefore, the sets of eigenvectors of $\tL^T\tL$, $\tL$, and $\tL^T$ are the same. 
Denote $V \in \mathbb{R}^{N\times N}$ the unitary matrix whose columns are eigenvectors of $\tL$ and $\Lambda \triangleq \di\{0,\lambda_{2},\ldots,\lambda_{N}\}$. 
%the diagonal matrix whose diagonal elements are eigenvalues of $\tL$. 
Consequently, we have $\tL=V\Lambda V^{*}$, $\tL^T=V\Lambda^{*}V^{*}$, $\tL^T\tL=V\Lambda^{*}\Lambda V^{*}$, and $P_{1}=V\di\{0,1,\ldots,1\}V^{*}$.  

Let $\lambda_{i,r}$ and $\lambda_{i,\ell}$ be the real and imaginary parts of $\lambda_{i},i=2,\ldots,N$, respectively.  
Then applying a congruence transformation with $V^{*}\otimes I_{n}$ to (\ref{Vdot-5}) gives us  
\begin{align}
	\label{Vdot-8}
	& XA^T+AX+\epsilon X %\nonumber \\
	 +\frac{1}{2}(\lambda_{i,r}+j\lambda_{i,\ell})B(\Delta_{1}+\Delta_{2})Y \nonumber \\
	& +B\Gamma^{-1}B^T+\frac{1}{2}(\lambda_{i,r}-j\lambda_{i,\ell})Y^T(\Delta_{1}+\Delta_{2})B^T \nonumber \\
	& +\frac{1}{4}(\lambda_{i,r}^{2}+\lambda_{i,\ell}^{2})Y^T\Gamma(\Delta_{1}-\Delta_{2})^2Y \preceq 0.
\end{align}
Furthermore, we have 
$	\lambda_{i,r} = 1-\cos\frac{2\pi(i-1)}{N}$ and $\lambda_{i,\ell} = -\sin\frac{2\pi(i-1)}{N}$ 
%\begin{align*}
	%\lambda_{i,r} = 1-\cos\frac{2\pi(i-1)}{N}, ~
	%\lambda_{i,\ell} = -\sin\frac{2\pi(i-1)}{N},
%\end{align*}
since $\tL$ is a circulant matrix with the following form
\begin{equation*}
	%\label{eq:}
	\tL = \begin{bmatrix} 1 & -1 & 0 & \cdots & 0 & 0 \\ 0 & 1 & -1 & \cdots & 0 & 0 \\ \vdots & \vdots & \ddots & & \vdots & \vdots \\ 0 & 0 & 0 & \cdots & 1 & -1 \\ -1 & 0 & 0 & \cdots & 0 & 1 
	      \end{bmatrix}.  
\end{equation*}
Denote 
$\mathbb{X}_{i} \triangleq XA^T+AX+\epsilon X+B\Gamma^{-1}B^T+\frac{1}{2}\lambda_{i,r}\mathrm{sym}(B(\Delta_{1}+\Delta_{2})Y)+\frac{1}{4}(\lambda_{i,r}^{2}+\lambda_{i,\ell}^{2})Y^T
\Gamma(\Delta_{1}-\Delta_{2})^2Y$. 
Subsequently, (\ref{Vdot-8}) is equivalent to 
$\begin{bmatrix} \mathbb{X}_{i} & -\mathbb{Y}_{i} \\ \mathbb{Y}_{i} & \mathbb{X}_{i} \end{bmatrix} \preceq 0.$ 
%\begin{align}
	%\label{Vdot-9}
	%\begin{bmatrix} \mathbb{X} & -\mathbb{Y} \\ \mathbb{Y} & \mathbb{X} \end{bmatrix} \preceq 0.
%\end{align}
On the other hand, using Schur complement for $\mathbb{X}_{i}$, we obtain
\begin{align}
	\label{Vdot-10}
	& \mathbb{X}_{i} \preceq 0  \nonumber \\
	\Leftrightarrow & \begin{bmatrix} \mathrm{sym}(AX+\frac{1}{2}B\tilde{\Delta}Y)+\epsilon X & BZ+Y^T\hat{\Delta}_{i} \\ Z^TB^T+\hat{\Delta}_{i}Y & -Z \end{bmatrix} \preceq 0,
\end{align}
where $Z \triangleq \Gamma^{-1}$.   
As a result, we obtain (\ref{input-LMI-d}). %(\ref{Vdot-9}) is equivalent to (\ref{input-LMI-d}). 

\end{pf}

\section{Numerical Example}
\label{app}

Consider a building temperature control problem where the target is to make the temperatures of different rooms be identical by allowing the exchange of their temperatures through a communication network. 
For simplicity, the dynamics of each room can be described by a first-order transfer function 
$\frac{a}{Ts+1},a>0,T>0,$ 
where the delays of the heating or cooling processes are ignored. 
Consequently, each room is equipped with an integrator controller for the consensusability, i.e., the model of each agent is $\frac{a}{s(Ts+1)}$.  
%then the state-space model of each room is represented by (\ref{agent}) with
%%$$G_{1}(s)=\frac{a}{s(Ts+1)}$$ with the associated   
%$$A_{1}=\begin{bmatrix} 0 & 1 \\ 0 & -\frac{1}{T} \end{bmatrix}, B_{1}=\begin{bmatrix} 0 \\ \frac{a}{T} \end{bmatrix}.$$

%--------------------------------------------------------
\subsection{Consensus under Input Constraints}

Then we illustrate the input-constrained consensus design in Theorem \ref{input-consensus-thm} and Theorem \ref{input-consensus-thm-1} 
 with a simulation of which $a=10,T=50$, $N=3$, and agents' inputs are bounded in $[-0.2,0.1]$. 

\begin{figure}[ht!]
		\centering
		\includegraphics[scale=0.4]{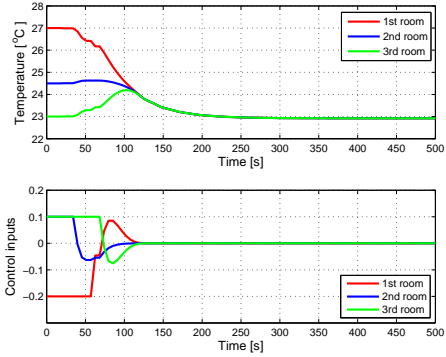}
		\caption{Consensus of temperatures under input constraint with {\it all-to-all undirected} communication structure.}
		\label{tempt-input-cnt-aut}
	\end{figure}

In the first scenario, the communication structure among rooms is undirected and all-to-all, and hence the eigenvalues of Laplacian matrix $\tL$ are $\{0,3,3\}$. 
Using $\epsilon=0.1$ and solving the LMI problem (\ref{input-LMI}) using {\tt CVX} \cite{Boyd:2015}, we obtain $\tK=[-0.1579,-1.0934]$.
Then the simulation results are shown in Figure \ref{tempt-input-cnt-aut}, which reveal that the temperatures of all rooms reach a consensus while the control inputs satisfy the given constraints. 
	
\begin{figure}[ht!]
		\centering
		\includegraphics[scale=0.4]{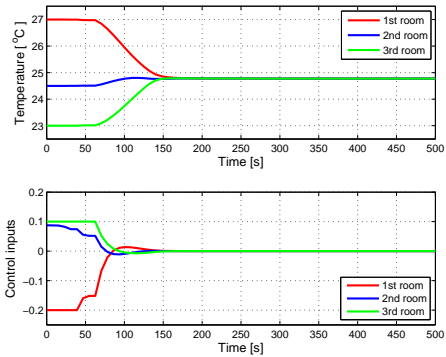}
		\caption{Consensus of temperatures under input constraint with an {\it undirected} communication structure.}
		\label{tempt-input-cnt-aut-1}
	\end{figure}

Consequently, we verify the agents' responses with a different undirected structure described by $2 \leftrightarrow 1 \leftrightarrow 3$, and the eigenvalues	of the associated Laplacian matrix 
are $\{0,1,3\}$.
Resolving the LMI problem (\ref{input-LMI}) using {\tt CVX} \cite{Boyd:2015} gives us $\tK=[-0.0350,-0.5946]$. 
Accordingly, the agents also reach consensus as seen in Figure \ref{tempt-input-cnt-aut-1}. However, the consensus value and the responses of agents are distinct with the above case of all-to-all 
undirected communication topology.

\begin{figure}[ht!]
		\centering
		\includegraphics[scale=0.4]{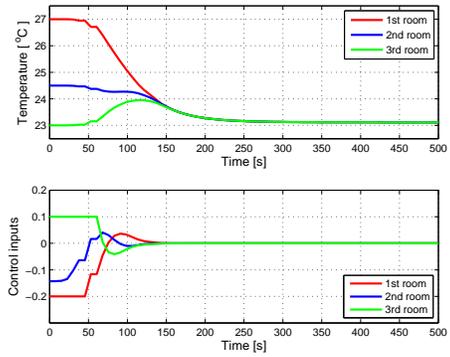}
		\caption{Consensus of temperatures under input constraint with {\it directed cyclic} communication structure.}
		\label{tempt-input-cnt-d}
	\end{figure}

Next, we consider the scenario where the interconnections among agents is described by a directed cyclic graph. 
In this case, we obtain $\tK=[-0.0953,-1.1379]$ by solving the LMI problem (\ref{input-LMI-d}) using {\tt CVX} \cite{Boyd:2015}. 
Then the simulation results are displayed in Figure \ref{tempt-input-cnt-d}. 
We can observe that the temperature of all rooms reach consensus despite the presence of the bounded input constraint and the directed communication topology. 
Furthermore, the consensus value, consensus speed, and the transient responses of rooms' temperatures are different from the previous cases of undirected communication structure.

%--------------------------------------------------------
\subsection{Consensus under Input Uncertainties }

\begin{figure}[ht!]
		\centering
		\includegraphics[scale=0.4]{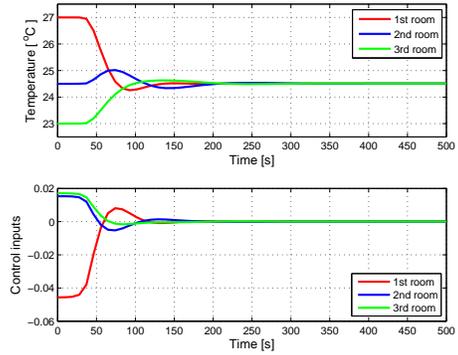}
		\caption{Consensus of temperatures under uncertain input with an {\it undirected} communication structure.}
		\label{tempt-input-unc-aut}
	\end{figure}

In this section, we assume that the control inputs of agents contain some uncertainties which may be due to the uncertain communication links. 
More specifically, $\delta_{k,1}$ and $\delta_{k,2}$ are assumed to be $0.7$ and $1.2$, respectively. 
This means the inputs of agents are multiplied with uncertain parameters $K_{1},K_{2},K_{3} \in [0.7,1.2]$. 
Then we utilize the same undirected structure described by $2 \leftrightarrow 1 \leftrightarrow 3$ as in the previous section and 
solve the LMI problem (\ref{input-LMI}) using {\tt CVX} \cite{Boyd:2015} to obtain $\tK=[-0.0059, -0.0760]$. 
Subsequently, we randomly generate the uncertainties on the inputs of agents and simulate the whole MAS to see the agents' responses. 
We observe that the consensus among agents is achieved for all uncertainties of the agents' inputs in the interval $[0.7,1.2]$.  
Results for an example with $K_{1}=1.1797$, $K_{2}=1.0279$, $K_{3}=0.7179$ are shown in Figure \ref{tempt-input-unc-aut}.

\begin{figure}[ht!]
		\centering
		\includegraphics[scale=0.4]{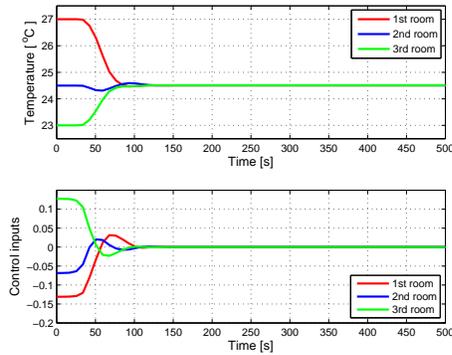}
		\caption{Consensus of temperatures under uncertain input with {\it directed cyclic} communication structure.}
		\label{tempt-input-unc-d}
	\end{figure}

Next, we investigate another context where the agents are interconnected through a directed cyclic graph. 
In this situation, solving the LMI problem (\ref{input-LMI-d}) using {\tt CVX} \cite{Boyd:2015} gives us $\tK=[-0.0444, -0.5154]$. 
Then we also see that the consensus of agents is obtained for any uncertainty of each agent's input in $[0.7,1.2]$. 
The simulation results with the same uncertain parameters $K_{1},K_{2},K_{3}$ as above are displayed in Figure \ref{tempt-input-unc-d} for comparison.

Overall, we can conclude that the consensus of agents under input constraints or uncertainties depends on the interconnection structure among agents. 
This is obvious since the communication topology affects to the eigenspectrum of the Laplacian matrix $\tL$ which directly influences the 
solutions of the LMI problems (\ref{input-LMI}) and (\ref{input-LMI-d}) and hence the consensus controller gain $\tK$. 
It is apparently different from the circumstance of consensus without any constraints or uncertainties since the consensus value is the average of initial 
conditions of agents interconnected by a connected undirected graph regardless of its structure.

%In another situation, we demonstrate the state-constrained consensus synthesis in Theorem \ref{state-consensus-thm} where agents' states are bounded by the intervals 
%$[-1,26]; [-1.5,1.5]$. 
%Choosing $\epsilon=0.1$ and solving the LMI problem (\ref{state-LMI}) gives us
%$K=\begin{bmatrix} -54.2360 & -289.3119 \end{bmatrix}.$ 	
%Consequently, the simulation results are displayed in Figure \ref{tempt-state-cnt}. 
%We can see that the temperatures of all rooms converge to a same value. 
	%\begin{figure}[ht!]
		%\centering
		%\includegraphics[scale=0.4]{tempt-state-cnt}
		%\caption{Consensus of temperatures with state constraint.}
		%\label{tempt-state-cnt}
	%\end{figure}	

%%%%%%%%%%%%%%%%%%%%%%%%%%%%%%%%%%%%%%%%%%%%%%%%%%%%%%%%%%%%%%%%%%%%%%%%%%%%%%%%	
\section{CONCLUSION}
\label{sum}

A new approach has been proposed in this paper to analyze and synthesize robust consensus controllers for linear leaderless MASs subjected to input constraints or uncertainties. 
The remarkable features of this approach are as follows. 
First, it is available for leaderless MASs with general linear dynamics of agents unlike the existing results for special cases of single integrator or single-input agents. 
Second, the robust consensus design under sector-bounded input constraints or uncertainties is derived in the form of a distributed low-dimension LMI problem which can be efficiently solved 
by off-the-shelf optimization software. 
Third, the proposed approach can deal with for both undirected and a special class of directed networks.  

The next researches would study more general classes of directed networks and take into account other practical issues such as time delays, disturbances, etc., together with 
the considered constraints and uncertainties.

%%%%%%%%%%%%%%%%%%%%%%%%%%%%%%%%%%%%%%%%%%%%%%%%%%%%%%%%%%%%%%%%%%%%%%%%%%%%%%%%
\section*{ACKNOWLEDGMENT}

This research is partially supported by Hitech Research Center, projects for private universities, supplied from the Ministry of Education, Culture, Sports, Science and Technology, Japan.

%%%%%%%%%%%%%%%%%%%%%%%%%%%%%%%%%%%%%%%%%%%%%%%%%%%%%%%%%%%%%%%%%%%%%%%%%%%%%%%%%%%%%%%%%%%%%%%%%%%%%%%%%%%%%
%%% The Appendices part is started with the command \appendix;
%\appendix
%
%%=========================================================================================
%\section{Proof of Theorem \ref{sync-n-ptb}}
%\label{A1}

%=======================================================================================================

\bibliographystyle{plain}        % Include this if you use bibtex 
\bibliography{References}           % and a bib file to produce the 
\end{document}